\documentclass[11pt, a4paper,]{article}
\usepackage{ amsfonts,amsmath}

\newtheorem{prop}{Proposition}[section]
\newtheorem{coro}{Corollary}[section]

\newcommand{\na}{\stackrel{\star}\nabla}
\newcommand{\sn}{\stackrel{\star}{A}_{\xi}}
\newtheorem{defi}{Definition}[section]
\newcommand{\ov}{\overline}
\newtheorem{theo}{Theorem}[section]
\newtheorem{lem}{Lemma}[section]
\newtheorem{rem}{Remark}[section]

\title{Einstein-Weyl structures on lightike hypersurfaces }
\author{C. Atindogbe\footnote{atindogb@iecn.u-nancy.fr,~  Permanent adress: Institut de mathematiques et de sciences Physiques (IMSP), Universit\'e d'Abomey-Calavi(UAC), Benin), 01 BP 613 Porto-Novo, Benin. Email: atincyr@imsp-uac.org} \qquad \qquad L. Berard-Bergery\footnote{berard@iecn.u-nancy.fr}\cr Institut Elie Cartan, Universit\'e Henri Poincar\'e, Nancy I, B.P. 239\cr 54506 Vand\oe uvre-l\`es Nancy Cedex, France}
\date{.}
\begin{document}
\maketitle

\begin{abstract}
\noindent
We study Weyl structures on lightlikes hypersurfaces endowed with a conformal structure of certain type and specific screen distribution: the Weyl screen structures. We investigate various differential geometric properties of Einstein-Weyl screen structures on  lightlike hypersurfaces and show that, for ambiant Lorentzian space~$\mathbb{R}^{n+2}_{1}$ and a totally umbilical screen foliation, there is a strong interplay with the induced (Riemannian) Weyl-structure on the leaves. Finally, we establish necessary and sufficient conditions for a Weyl structure defined by the $1-$form of an almost contact structure given by an additional complex structure  in case of an ambiant Kaehler manifold to be closed.
\end{abstract}

\noindent
{\bf{Key words:~}}Lightlike hypersurface, screen distribution,  Einstein-Weyl structure.

\vspace{0.5cm}
\noindent
{\bf{MSC subject classification (2000):~}} 53C50, 53C05, 53C25.             

\section{Introduction }
\label{intro}
Pseudo-Riemannian manifolds $(M,g)$ with $dim M = n > 4$ and $sgn(g)= (n-1,1)$ are natural generalizations of ($4$-dimensional model) spacetime of general relativity. Lightlike hypersurfaces in $(M,g)$ are models of different types of horizons separating domains of $(M,g)$ with different physical properties.

As it is well known, contrary to timelike and spacelike hypersurfaces, the geometry of  lightlike hypersurfaces is different and rather difficult since the normal bundle and the tangent bundle have non-zero intersection.

Being lighlike manifold is invariant under conformal change of the metric, along with many geometric objects. Thus, it is more suitable to study geometry of lightlike hypersurfaces within the framework of conformal class of degenerate metrics.

In a Riemannian setting, manifolds $M^{n}$  with conformal structure $[g]$ and torsion-free connection $D$, such that parallel translation induces conformal transformations, are called Weyl manifolds. They are said to be Einstein-Weyl if the symmetric trace-free  part of the Ricci tensor of the (Weyl) connection $D$ vanishes. If $D$ is locally the Levi-Civita connection of a compatible metric in $[g]$, the structure is said to be closed, and the ($D$-compatible) metric is locally Einstein \cite{Gaud,Pedswan, Iva}. 

In \cite{DB}, Duggal and Bejancu introduced a main tool in studying the geometry of a lightlike hypersurface: the screen distributions. The latter is used to constrtuct a lightlike transversal vector bundle which is non-intersecting to the lightlike tangent bundle. It is now well-known  that a suitable choice of screen distribution has produced substantial result in lightlike geometry\cite{DB,AD}. Based on this, we brief in section~\ref{prelim} basic informations on normalizations, induced geometric objects \cite{DB} and pseudo-inversion of degenerate metrics \cite{ATE}.In section~\ref{weyl}, we define \emph{Weyl screen structure} (Definition~\ref{weyl1}),  and prove a result on model space of Weyl screen structures on the (conformal) lightlike hypersurface. Therefater, we study and relate curvature and Ricci tensors of the Weyl connection, along with its scalar curvature to their respective analogous for a given representative element in the conformal class.  In section~\ref{einsweyl}, we consider Einstein-Weyl screen structures and establish a necessary and sufficient condition for a Weyl screen structure to be Einstein-Weyl. Section~\ref{totumb}  is devoted to a special case of total umbilicity of the screen foliation involved in Definition~\ref{weyl1}. Also, in ambiant Lorentzian case, we prove that there is a strong interplay between Einstein-Weyl screen  structures on the conformal lightlike hypersurface and the (induced) one on the (Riemannian) screen foliation. Section~\ref{kahler}  deals with lightlike real hypersurfaces of Kahler manifolds.

\section{Preliminaries on Lightlike hypersurfaces}
\label{prelim}
Let $M$ be a hypersurface  of an $(n+2)-$dimensional pseudo-Riemannian manifold $(\ov{M},\ov{g})$~of index $0 < \nu < n+2 $. In the classical theory of nondegenerate hypersurfaces, the normal bundle has trivial intersection $\{0\}$ with the tangent one and plays an important role in the introduction of main geometric objects. In case of lightlike (degenerate, null) hypersurfaces, the situation is totally different. The normal bundle $TM^{\perp}$ is a rank-one distribution over $M$: $TM^{\perp}\subset TM$ and then coincide with the so called \emph{radical distribution} $RadTM = TM \cap TM^{\perp}$.  Hence,the induced metric tensor field $g$ is degenerate and has rank $n$. The following characterisation is proved in \cite{DB}.   

\begin{prop}  
\label{prop1}
Let $(M,g)$ be a hypersurface of an $(n+2)-$dimensional pseudo-Riemannian manifold $(\ov{M},\ov{g})$. Then the following assertions are equivalent.
\begin{itemize}
\item[(i)]
$M$ is a lightlike hypersurface of $\ov{M}$.
\item[(ii)]
$g$ has constant rank $n$ on $M$.
\item[(iii)]
$TM^{\perp}= \cup_{x\in M}T_xM^\perp$ is a distribution on $M$.
\end{itemize}
\end{prop}
A complementary  bundle of $TM^\perp$ in $TM$ is a rank $n$ nondegenerate distribution over $M$. It is called a \emph{screen distribution} on $M$  and is often denoted by $S(TM)$. A lightlike hypersurface endowed with a specific screen distribution is denoted by the triple $(M,g,S(TM))$. As $TM^\perp$ lies in the tangent bundle, the following result has an important role in studyng the geometry of a lightlike hypersurface.

\begin{prop}(\cite{DB})
\label{theo1}
Let $(M,g,S(TM))$ be a lightlike hypersurface of $(\ov{M},\ov{g})$ with a given screen distribution $S(TM)$. Then there exists a unique rank $1$ vector subbundle $tr(TM)$ of $T\ov{M}|_{M}$, such that for any non-zero section $\xi$ of $TM^\perp$ on a coordinate neighbourhood ${\mathcal{U}}\subset M$, there exists a unique section $N$ of $tr(TM)$ on ${\mathcal{U}}$ satisfyng 
\begin{equation}
\label{eq1}
\ov{g}(N,\xi)= 1
\end{equation}
and
\begin{equation}
\label{eq2}
\ov{g}(N,N)= \ov{g}(N,W)= 0, \quad \quad \forall~ W \in \Gamma(ST|_{\mathcal{U}}).
\end{equation}
\end{prop}
Here and in the sequel we denote by $\Gamma(E)$ the ${\mathcal{F}}(M)-$module of smooth sections of a vector bundle $E$ over $M$, ${\mathcal{F}}(M)$ being the algebra of smooth functions on $M$. Also, by $\perp$ and $\oplus$ we denote the orthogonal and non-orthogonal direct sum of two vector bundles.  By proposition~\ref{theo1} we may write down the following decompositions.

\begin{equation}
\label{eq3}
TM=S(TM) \perp TM^\perp,
\end{equation}

\begin{equation}
\label{eq4}
T\ov{M}|_{M} = TM \oplus tr(TM)
\end{equation}

and

\begin{equation}
\label{eq4bis}
T\ov{M}|_{M}= S(TM) \perp (TM^\perp \oplus tr(TM))
\end{equation}

As it is well known, we have the following:
\begin{defi}
\label{induced}
Let $(M,g,S(TM))$ be a lightlike hypersurface of $(\ov{M},\ov{g})$ with a given screen distribution $S(TM)$. The induced connection, say $\nabla$, on $M$ is defined by
\begin{equation}
\label{eq48}
\nabla_{X}Y = Q(\ov{\nabla}_XY),
\end{equation} 
where $\ov{\nabla}$ denotes  the Levi-civita connection on $(\ov{M},\ov{g})$ and $Q$ is the projection on $TM$ with respect to the decomposition $(\ref{eq4})$.
\end{defi}

\begin{rem}
\label{rem1}
Notice that the induced connection $\nabla$ on $M$ depends on both $g$ and the specific given screen distribution $S(TM)$ on $M$.
\end{rem}
By respective projections $Q$ and $I-Q$, we have Gauss an Weingarten formulae in the form

\begin{equation}
\label{eq5}
\ov{\nabla}_XY= \nabla_XY + h(X,Y)\qquad \forall X,Y ~ \in \Gamma(TM),
\end{equation}

\begin{equation}
\label{eq6}
\ov{\nabla}_XV= -A_VX + \nabla_X^tV \qquad \forall X ~ \in \Gamma(TM),\quad \forall ~V \in \Gamma(tr(TM)).
\end{equation}
Here, $\nabla_XY$ and $A_{V}X$ belong to $\Gamma(TM)$. Hence

$\bullet$  $h$ is a $\Gamma(tr(TM))$-valued symmetric ${\mathcal{F}}(M)$-bilinear form on $\Gamma(TM)$,

$\bullet$ $A_{V}$ is an  ${\mathcal{F}}(M)$-linear operator on $\Gamma(TM)$, and 

$\bullet$ $\nabla^t$ is a linear connection on the lightlike transversal vector bundle $tr(TM)$. 

Let $P$  denote the projection morphism of $\Gamma(TM)$ on $\Gamma(S(TM))$ with respect to the decomposition (\ref{eq3}). We have 

\begin{equation}
\label{eq9}
\nabla_XPY= \na_XPY + h^{*}(X,PY)\qquad \forall X,Y ~ \in \Gamma(TM),
\end{equation}

\begin{equation}
\label{eq10}
\nabla_X U= -\stackrel{\star}{A}_{U}X + \nabla^{*t}_XU \qquad \forall X ~ \in \Gamma(TM),\quad \forall ~U \in \Gamma(TM^\perp).
\end{equation}

Here $\na_XPY$ and $\stackrel{\star}{A}_{U}X$ belong to $\Gamma(S(TM))$,  $\na$ and $\nabla^{*t}$are linear connection on $S(TM)$ and $TM^\perp$, respectively. Hence

$\bullet$  $h^{*}$ is a $\Gamma(TM^\perp)$-valued  ${\mathcal{F}}(M)$-bilinear form on $\Gamma(TM)\times \Gamma(S(TM))$, and 

$\bullet$ $\stackrel{\star}{A}_{U}$ is a $\Gamma(S(TM))$-valued ${\mathcal{F}}(M)$-linear operator on $\Gamma(TM)$.
 
\noindent
They are the second fundamental form  and the shape operator of the screen distribution, respectively. 

Equivalently, consider a normalizing pair $\{\xi, N\}$ as in the proposition~\ref{theo1}. Then, $(\ref{eq5})$ and $(\ref{eq6})$ take the form

\begin{equation}
\label{eq7}
\ov{\nabla}_XY= \nabla_XY + B(X,Y)N \qquad \forall X,Y ~ \in \Gamma(TM|_{\mathcal{U}}),
\end{equation}
and
\begin{equation}
\label{eq8}
\ov{\nabla}_XN= -A_{N}X+ \tau(X)N \qquad \forall X ~ \in 
\Gamma(TM|_{\mathcal{U}}),
\end{equation}
where we put locally on ${\mathcal{U}}$,
\begin{equation}
\label{eq13}
B(X,Y) = \ov{g}(h(X,Y),\xi)
\end{equation}

\begin{equation}
\label{eq14}
\tau(X) = \ov{g}(\nabla^{t}_XN,\xi)
\end{equation}

\noindent
It is important to stress the fact that the local second fundamental form $B$ in $(\ref{eq13})$ does not depend on the choice of the screen distribution.

We also define (locally) on ${\mathcal{U}}$ the following:
\begin{equation}
\label{eq11}
C(X,PY) = \ov{g}(h^{*}(X,PY),N),
\end{equation}

\begin{equation}
\label{eq12}
\varphi(X) = - \ov{g}(\nabla^{\star t}_X\xi ,N).
\end{equation}

Thus, one has for $~X\in~\Gamma(TM)$
\begin{equation}
\label{eq15}
\nabla_XPY= \na_XPY + C(X,PY)\xi
\end{equation}

\begin{equation}
\label{eq16}
\nabla_X\xi = -\sn X + \varphi(X)\xi 
\end{equation}
It is straighforward to verify that for $~X,Y \in~\Gamma(TM)$
\begin{equation}
\label{eq17}
B(X,\xi) = 0 
\end{equation}

\begin{equation}
\label{eq18}
B(X,Y) = g(\sn X,Y)
\end{equation}

\begin{equation}
\label{eq19}
\sn\xi = 0
\end{equation}

The linear connection $\na$ from (\ref{eq9})is a metric connection on $S(TM)$ and we have for all tangent vector fields $X$, $Y$ and $Z$ in $TM$
\begin{equation}
\label{eq20}
\left(\nabla_{X}g\right)(Y,Z)~=~B(X,Y)\eta(Z) + B(X,Z)\eta(Y).
\end{equation}
with
\begin{equation}
\label{eq49}
\eta(\cdot) = \ov{g}(N,\cdot).
\end{equation}

The induced connection $\nabla$ is torsion-free, but not necessarily $g$-metric. Also, on the geodesibility of $M$ the following is known.

\begin{theo}(\cite[p.88]{DB})
\label{theo2}
Let $(M,g,S(TM))$ be a lightlike hypersurface of a pseudo-Riemannian manifold $(\ov{M},\ov{g} )$. Then the following assertions are equivalent:
\begin{itemize}
\item[(i)] $M$ is totally geodesic.
\item[(ii)]$h$ (or equivalently $B$) vanishes identically on $M$.
\item[(iii)] $\stackrel{\star}{A}_{U}$ vanishes identically on $M$, for any $U~\in \Gamma(TM^\perp)$
\item[(iv)]The  connection $\nabla$ induced by $\ov{\nabla}$ on $M$ is torsion-free and metric.
\item[(v)] $TM^\perp$ is a parallel distribution with respect to $\nabla$.
\item[(vi)]$TM^\perp$ is a Killing distribution on $M$.
\end{itemize}
\end{theo}
It turns out that if $(M,g)$ is not totally geodesic, there is no connection that is, at the same time, torsion-free and $g$-metric.

\subsection{Pseudo-inversion of degenerate metrics}
A large class of differential operators in differential geometry is intrinsically defined by means of the dual metric $g^{*}$ on the dual bundle $\Gamma(T^{*}M)$~of $1$-forms on $M$. If the metric $g$ is nondegenerate, the tensor field $g^{*}$ is nothing but the inverse of $g$. We brief here construction of some of these operators in case the metric $g$ is degenerate and refer the reader to \cite{ATE} for more details.

Let $(M,g,S(TM))$ be a lightlike hypersurface and $\{\xi,N  \}$ be a pair of  (null-) vectors given by the normalizing Proposition~ \ref{theo1}. Consider on $M$ the one-form defined by
\begin{equation} 
\label{eq20bis}
\eta(\cdot)~=~ \ov{g}(~N~,~\cdot~)
\end{equation}
For all $X\in \Gamma(TM)$,
$$X = PX + \eta(X)\xi $$ 
and $\eta(X)= 0 $ if and only if $X\in \Gamma(S(TM))$. Now, we define $\flat$ by 

\begin{eqnarray*}
\flat: \Gamma(TM)  &  \longrightarrow & \Gamma(T^{*}M)\cr
                   &                  &               \cr 
                X  &  \longmapsto         &  X^{\flat}
\end{eqnarray*}
such that 
\begin{equation} 
\label{eq21}
X^{\flat} = g(~X~, ~\cdot~) + \eta(X)\eta(~\cdot~)
\end{equation}

Clearly, such a $\flat$ is an  isomorphism of $\Gamma(TM)$ onto $\Gamma(T^{*}M)$, and generalize the usual nondegenerate theory. In the latter case,  $\Gamma(S(TM))$ coincide with 
$\Gamma(TM)$, and as a consequence the $1-$form $\eta$ vanishes identically and the projection morphism P becomes the identity map on $\Gamma(TM)$. We let $\sharp$ denote the inverse of the isomorphism $\flat$ given by (\ref{eq21}). For $X\in \Gamma(TM)$ (resp. $\omega \in T^{*}M $), $X^{\flat}$ (resp. $\omega^{\sharp}$) is called the dual $1-$form of $X$ (resp. the dual vector field of $\omega$) with respect to the degenerate metric $g$. It follows (\ref{eq21}) that if $\omega$ is a $1$-form on $M$, we have for $X\in \Gamma(TM)$

\begin{equation}
\label{eq33}
\omega(X)~=~g(\omega^{\sharp},X) ~+~\omega(\xi)\eta(X)
\end{equation}

Now we introduce the so called associate non degenerate metric $\tilde{g}$ to the degenerate metric $g$ as follows. For $X,Y \in \Gamma(TM)$, define $\tilde{g}$ by 

\begin{equation} 
\label{eq22}
\tilde{g}(X,Y)~=~X^{\flat}(Y) 
\end{equation}
Clearly,$ \tilde{g} $ defines a non degenerate metric on $M$ and play an important role in defining the usual differential operators \emph{gradient}, \emph{divergence}, \emph{laplacian} with respect to degenerate metric $g$ on lightlike hypersurfaces. Also, obseve that  $\tilde{g}$ coincides with $g$ if the latter is not degenerate. The $(0,2)$ tensor field $g^{[~\cdot~,~\cdot~]}$, inverse of $\tilde{g}$ is called \emph{the pseudo-inverse of $g$}. Finally, we state the following result (\cite{ATE}).
\begin{prop}
\label{prop2}
Let $(M,g,S(TM))$ be a lightlike hypersurface of a pseudo-Riemannian $(n+2)$-dimensional  manifold $(\ov{M},\ov{g} )$.We have 
\begin{itemize}
\item[(i)] for any smooth function $f:{\mathcal{U}}\subset M \rightarrow \mathbb{R}$,
\begin{equation}
\label{eq23}
grad^{g}f~=~ g^{[\alpha \beta]}f_{\alpha}\partial_\beta 
\end{equation}
where $f_{\alpha}= \frac{\partial f}{\partial x^{\alpha}}$, $\partial_{\beta}= \frac{\partial}{\partial x^{\beta}}$,~~$\alpha,~\beta= 0,\cdots, n $;
\item[(ii)] For any vector field $X$ on ${\mathcal{U}}\subset M$,
\begin{equation}
\label{eq24}
div^{g}X = \sum_{\alpha,\beta}g^{[\alpha,~\beta]}\tilde{g}(\nabla_{\partial_{\alpha}}X,\partial_{\beta})
\end{equation}
\item[(iii)] for smooth function $f:{\mathcal{U}}\subset M \rightarrow \mathbb{R}$
\begin{equation}
\label{eq25}
\Delta^{g}f~=~ \sum_{\alpha \beta}g^{[\alpha,~\beta]}\tilde{g}(\nabla_{\partial_{\alpha}}grad^{g}f,\partial_{\beta})
\end{equation}
\end{itemize}
where $\{\partial_{0}:=\xi, \partial_{1},\cdots,\partial_{n} \}$ is any quasiorthonormal frame field on $M$ adapted to the decomposition (\ref{eq3}).
\end{prop}

In index free notation, (\ref{eq23}) can be written in the form ~$\tilde{g}(\nabla^{g}f,X)~=~df(X)$ which defines the gradient of the scalar function $f$ with respect to the degenerate metric $g$. With nondegenerate $g$, one has $\tilde{g} = g$ so that $(i)-(iii)$ generalize the usual known formulae to the degenerate set up. 

From now on, unless otherwise stated, the ambiant manifold $(\ov{M},\ov{g})$ has a Lorentzian signature so that all lighlike hypersurfaces considered are of signature $(0,n)$. In particular, it follows that any screen distribution is Riemannian. 

As it is well known (theorem~\ref{theo2}), only totally geodesic lightlike hypersurfaces do have their induced connection metric and torsion-free. In the next section and the remainder of the text, only such lightlike hypersurfaces will be in consideration. Also, being lightlike is invariant under conformal change of the metric. Athough, for $(M,g_{0})$ totally geodesic, not all metrics in the conformal class of $g_{0}$  guarantee this geometric condition. In this respect, we consider appropriate conformal structure on a given totally geodesic $(M,g_{0})$.

\section{Weyl screen structures}
\label{weyl}
Let $(M,g_{0})$ be a totally geodesic hypersurface in a $(n+2)-$dimensional pseudo-Riemannian manifold $(\ov{M},\ov{g})$. Consider on $M$ conformal metrics of the form $g= e^{-2f}g_{0}$ with $X(f)=0$ for $X\in TM^{\perp} = span\{\xi \}$ i.e $f$ is constant on $\xi-orbits$. These metrics endow $M$ with a special conformal structure we denote by $c=[g_{0}]_{0}$. For each metric $g \in c$, $(M,g)$ is also totally geodesic, and there exists a $g-$compatible torsion-free connection $\nabla^{g}$. Throughout the text, $M$, endowed with this conformal  structure is denoted $(M,c)$.

\begin{defi}
\label{preweyl}
A Weyl structure relative to $(M,c)$ is a symmetric linear connection $D$ on $M$ that preserves the structure. More precisely, $D$ satisfies

\begin{itemize}
\item[(i)] $D$ is torsion-free
\item[(ii)] For $g$ in the conformal class $c$, there exists a unique $1-$form $\theta$ on $M$ such that 
\begin{equation}
\label{eq26}
Dg = -2\theta \otimes g
\end{equation}
\end{itemize}
\end{defi}

\begin{rem}
\label{rem1bis}
Conditions $(i)$ and $(ii)$ in definition~\ref{preweyl} determine a Weyl structure modulo $S^{2}T^{*}M\otimes TM^{\perp}$
\end{rem}

\begin{lem}
\label{lem1}
The Kernel $TM^{\perp}~(=RadTM = Ker~g)$ of $g$ is parallel with respect to any Weyl structure $D$ on $(M,c)$. Furthermore, up to a renormalization, one can choose $\xi \in TM^{\perp}$ so that $D_\xi \xi =0$ and for any $g \in c$, there exists a torsion-free $g$-compatible linear connection $D^{g}$ with $D^{g}_\xi \xi =0$.
\end{lem}

\noindent
{\bf{Proof}}
Let $X,Y,Z \in \Gamma(TM)$ and $g\in c$. From (\ref{eq26}) we have
$$X\cdot g(Y,Z)-g(D_{X}Y,Z)-g(Y,D_{X}Z)= -2\theta(X)~g(Y,Z)  $$
Then for  $Z\in RadTM, $ one has $g(Y,D_{X}Z) = 0~ \forall ~Y\in \Gamma(TM)$. Thus  $D_{X}Z \in RadTM~\forall ~ Z\in RadTM$. Now let $\xi \in RadTM $, we have $D_\xi \xi = \psi(\xi)\xi$. If $\psi(\xi) = 0$ then there is nothing more to prove. Otherwise, choose on the null integral curve ${\mathcal{C}}$of $\xi$ a new parameter $t^{*}(t)$ such that 
$$\frac{d^{2}t^{*}}{dt^{2}}-\psi(\frac{d}{dt})\frac{dt^{*}}{dt} = 0$$
with $\frac{d}{dt}=:\xi$. Such a parameter always exists on ${\mathcal{C}}$ and one has $D_{\frac{d}{dt^{*}}}\frac{d}{dt^{*}}=0$. Now, let $g\in c$ and $D_{1}^{g}$ be a torsion-free $g$-compatible connection. Then, $0 = D_\xi \xi=D_{1\xi}^{g} \xi + S(\xi ,\xi)\xi$ where $S\in S^{2}T^{*}M$. If $D_{1\xi}^{g} \xi = 0$ then there is nothing more to prove. Otherwise, change $D_{1}^{g}$ in $D_{2}^{g}=D_{1}^{g} + S \otimes \xi$. Such a $D_{2}^{g}$ is a torsion-free linear $g$-compatible connection on $M$ and  
${D_{2}}_\xi^{g} \xi = 0$ and the proof is complete.$\Box$

\begin{rem}
\label{rem2}
From lemma~\ref{lem1} it follows that the element $S\in S^{2}T^{*}M$ modulo which the Weyl structure is determined satisfies $S(\xi,\xi)=0$ for a suitable choice of the torsion-free  $g$-compatible metric $D^{g}$ of $g$. The element $S\in S^{2}T^{*}M$ is entirely determined by the following. 
\end{rem}

\begin{defi}
\label{weyl1}
Let $(M,c,S(TM))$ be a totally geodesic lightlike hypersurface $(M,g_{0})$ endowed with the conformal structure $c=[g_{0}]_{0}$, and an integrable screen distribution $S(TM)$. A \emph{Weyl screen structure} $D$ relative to $(M,c,S(TM))$ is a Weyl structure for which $S(TM)$ is parallel, that is for all tangent vector fields $X$ and $Y$ in $TM$,
\begin{equation}
\label{eq27}
D_{X}PY ~\in \Gamma(S(TM))
\end{equation}
\end{defi}

{\bf{Note.~}} Throughout the text, we sometimes  consider the quadruplet $(M,c,D,S(TM))$ (as in Definition~\ref{weyl1}) as the  Weyl screen structure. Also, vector fields tangent to leaves of the refered screen distribution are called \emph{horizontal}.

\begin{lem} Let $D$ be a Weyl screen structure on $(M,c,S(TM))$.
\label{lem2}
Let $\Omega^{1}_{hor}(M)$ denote the space of horizontal $1-$form on $M$, that is $\omega \in \Omega^{1}_{hor}(M)$ if and only if $\omega (X)=0$ for all $X\in RadTM$.
\begin{itemize}
\item[(i)] For any  $g\in c$, ~ $\theta_{g} \in \Omega^{1}_{hor}(M)$.
\item[(ii)] For $g \in c$ there exists a unique $\theta_{g} \in \Omega^{1}_{hor}(M)$and a unique $S\in S^{2}T^{*}M$ such that for $X,Y \in \Gamma(TM)$,
\begin{equation}
\label{eq28}
D_{X}Y=D^{g}_{X}Y +\theta_{g}(X)Y +\theta_{g}(Y)X - g(X,Y)\theta_{g}^{\sharp_{g}} - S(X,Y)\xi
\end{equation}
where $\theta_{g}^{\sharp_{g}}$ is the dual of $\theta_{g}$ with  respect to the degenerate metric $g$ and the screen distribution $S(TM)$. Furthermore,

\begin{eqnarray}
\label{eq29}
S(X,Y)&=& \left\{\begin{array}{lcl} 0 & \mbox{if}& X,Y \in RadTM\cr
& & \cr
C(X,Y)+\eta(X)\theta_{g}(Y) &\mbox{if}& (X,Y)\in \cr & &\Gamma(TM)\times \Gamma(S(TM))
   \end{array}\right.
\end{eqnarray}
where $C$ denotes the second fundamental form of $S(TM)$ in $(M,g)$.
\end{itemize}
\end{lem}

\noindent
{\bf{Proof}}
Let $X\in RadTM$, $Y,Z \in \Gamma(TM)$. From (\ref{eq26}) and lemma~\ref{lem1} we have $L_{X}g_{0}(Y,Z)=-2~\theta_{g_{0}}(X)g_{0}(Y,Z)$. But $(M,g_{0})$ is totally geodesic and $L_{X}g_{0} = 0$~(theorem~\ref{theo2}). Thus,  $\theta_{g_{0}}(X)=0 ~\quad X\in RadTM$. For $g = e^{-2f}g_{0}\in c$,we have $\theta_{g}= \theta_{g_{0}}+ df$ with $df(X)=0~\forall~ X\in RadTM$. Thus, 
$\theta_{g}(X)= \theta_{g_{0}}(X)+ df(X)=0~\forall~X\in ~RadTM $ and $(i)$ is proved.

Now,let us write  for a choice of $g\in c$ ~and for all ~$X,Y\in \Gamma(TM)$

\begin{equation}
\label{eq31}
D_{X}Y~=~D^{g}_{X}Y~+~\tilde{\theta}_{X}Y
\end{equation}
where $D^{g}_{X}Y$ is the torsion-free $g$-compatible linear connection pointed out in lemma~\ref{lem1}. As $D$ and $D^{g}$ are torsion-free, one has 
\begin{equation}
\label{eq30}
\tilde{\theta}_{X}Y~=~\tilde{\theta}_{Y}X
\end{equation}

Taking into account (\ref{eq31}), (\ref{eq30}) and the $g$-compatibility of $D^{g}$ one has 

\begin{equation}
\label{eq32}
g(\tilde{\theta}_{X}Y,Z)~+~g(Y,\tilde{\theta}_{X}Z)~=~2\theta_{g}(X)g(Y,Z)
\end{equation}
By circular permutation in (\ref{eq32}) and taking into account (\ref{eq30}) one has
$$g(\tilde{\theta}_{X}Y,Z)~=~\theta_{g}(X)g(Y,Z)+\theta_{g}(Y)g(X,Z)-\theta_{g}(Z)g(X,Y)$$
As $\theta_{g}$ is horizontal (from (i)) its $g-$dual $\theta^{\sharp_{g}}_{g}$ is a horizontal vector field and from (\ref{eq33}) one can write ~$\theta_{g}(Z)~=~ g(Z,\theta^{\sharp_{g}}_{g})$. It follows that

$$\tilde{\theta}_{X}Y ~=~\theta_{g}(X)Y+\theta_{g}(Y)X - g(X,Y)\theta^{\sharp_{g}}_{g}-S(X,Y)\xi$$

for some $S\in S^{2}T^{*}M$. Also, from (\ref{eq15}) we have
$$D^{g}_{X}PY~=~\na^{g}_{X}PY ~+~C(X,PY)\xi $$
where $\na^{g}$ is the induced Levi-Civita connection by $D^{g}$ on the screen distribution  and $C$ the second fundamental form of the screen distribution in $(M,g)$. Thus

\begin{eqnarray}
\label{eq34}
D_{X}PY &=&\na^{g}_{X}PY +  \theta_{g}(X)PY + \theta_{g}(Y)PX -g(X,Y)\theta^{\sharp_{g}}_{g}\cr     & &+ \left[C(X,PY)+\eta(X)\theta_{g}(Y)-S(X,PY)\right]\xi
\end{eqnarray}
Observe that, since $\theta_{g}$ is a horizontal $1-$form, one has $\theta^{\sharp_{g}}_{g} \in \Gamma(S(TM))$. From $(iii)$ in definition~\ref{eq27}, $S(TM)$ is $D-$parallel if and only if the term in bracket vanihes identically on $M$. it follows that for 

\begin{equation}
\label{eq35}
X,Y \in \Gamma(TM), ~~S(X,PY)=C(X,PY)+\eta(X)\theta_{g}(Y)
\end{equation}
In particular

\begin{equation}
\label{eq36}
\noindent
\forall~Y \in \Gamma(TM), ~~S(\xi,PY)=S(PY,\xi)= C(\xi,PY)+\theta_{g}(Y)
\end{equation}
Finally, $S(\xi,\xi)=0$ follows remark~\ref{rem2} and the proof is complete.

\begin{rem}
\label{rem3}
\begin{itemize}
\item[(a)] From $S(\xi,\xi)=0$ and (\ref{eq36}) one can write
\begin{equation}
\label{eq37}
\forall ~Y \in \Gamma(TM), ~~S(\xi,PY)=S(PY,\xi)=: C(\xi,PY)+\theta_{g}(Y).
\end{equation}
\item[(b)] Clearly, for a given $g\in c$,~among all $g-$compatible torsion-free linear connections, there is only one which satisfies (\ref{eq28}). Thus, if we take our data for a Weyl screen structure on $(M,S(TM))$ to be $g \in c$ and the $1-$form $\theta_{g}$, $D= D^{g} + \tilde {\theta}$ is uniquely determined.
\end{itemize}
\end{rem}

The curvature tensor of the Weyl screen structure $D$ is defined by 
\begin{equation}
\label{eq41}
R^{D}(X,Y)~=~D_{[X,Y]}~-~[D_{X},D_{Y}]
\end{equation}
and we let $Ric^{D}$ denote the Ricci curvature of $D$. It is defined to be  the trace  of the map  $~Z~\mapsto~R^{D}(X,Z)Y$.~ For a representative $g\in c$ and a $g-$quasiorthonornal frame field $(X_{\alpha})_{\alpha}$ on $M$, 
\begin{equation}
\label{eq38}
Ric^{D}(X,Y)~=~g^{[\alpha \beta]}\tilde{g}(R^{D}(X,X_{\alpha})Y,X_{\beta})
\end{equation}

and clearly, the right hand side of (\ref{eq38}) does not change under conformal rescaling in $c$. The scalar curvature $Scal^{D}$ of $D$ is defined by 
$$Scal^{D}~=~tr_{c}(Ric^{D})$$
Observe that $Scal^{D}$ is not a function on $M$, but for a choice of a metric $g\in c$, it is   defined by

\begin{equation}
\label{eq39}
Scal^{D}_{g}~=~tr_{g}(Ric^{D})
\end{equation}

\begin{prop}
\label{curv}
Suppose $D~=~D^{g}~+~\tilde{\theta}$ where $g\in c$ and $\theta_{g}$ the $1-$form associated to the paire $\{D,g \}$. Then
\begin{eqnarray}
\label{eq40}
R^{D}(X,Y)&=&R^{g}(X,Y)-2d\theta_{g}(X,Y)id  \cr 
& & + \left(D^{g}_{Y}\theta_{g}^{\sharp_{g}}-\theta(Y)\theta_{g}^{\sharp_{g}} + \frac{1}{2}|\theta_{g}^{\sharp_{g}}|^{2}_{g}Y \right)\wedge X \cr
& & - \left(D^{g}_{X}\theta_{g}^{\sharp_{g}}-\theta(X)\theta_{g}^{\sharp_{g}} + \frac{1}{2}|\theta_{g}^{\sharp_{g}}|^{2}_{g}X \right)\wedge Y \cr 
& &\cr 
& & -
\left({\mathcal{K}}^{g}(X,Y)- {\mathcal{K}}^{g}(Y,X) \right)\xi 
\end{eqnarray}
with 
$${\mathcal{K}}^{g}(X,Y)= i_{Y}(D^{g}_{X}S) + S(Y,\theta_{g}^{\sharp_{g}})i_{X}g+S(Y,\xi)i_{X}S + \varphi_{g}(X)i_{Y}S$$
where ($(M,g)$ being totally geodesic) the $1-$form $\varphi_{g}$ is defined by $D^{g}_{X}\xi = \varphi_{g}(X)\xi$, and $X\wedge Y = g(X,\cdot)Y - g(Y,\cdot)X$.
\end{prop}

{\bf{Proof}} This is a standard computation using (\ref{eq28}) and the curvature formula (\ref{eq41}).$\Box$

The following lemma gives expression of ${\mathcal{K}}^{g}(X,Y)- {\mathcal{K}}^{g}(Y,X)$ for horizontal $X$ and $Y$ in terms of the second fundamental form $C$ of the screen distribution $S(TM)$. 
\begin{lem}
\label{matcalK}
For $X,Y\in \Gamma(S(TM))$, we have
\begin{eqnarray}
{\mathcal{K}}^{g}(X,Y)- {\mathcal{K}}^{g}(Y,X)&=&\eta(\bar{R}(X,Y)Z)\cr & &+ \left[g(X,Z)c(Y,\theta_{g}^{\sharp_{g}})- g(Y,Z)c(X,\theta_{g}^{\sharp_{g}}) \right] \cr & & + \left[C(X,Z)C(\xi,Y)- C(Y,Z)C(\xi,X) \right] \cr & &+
\left[\theta_{g}(Y)C(X,Z)-\theta_{g}(X)C(Y,Z)  \right]
\end{eqnarray}

where $\bar{R}$ is the ambiant Riemannian curvature of $(\ov{M},e^{-2\bar{f}}\ov{g})$, with $\bar{f}|_{M}=f$ and $C$ the second fundamental form of the screen distribution $S(TM)$.
\end{lem}

{\bf{Proof}~}This  is a direct use of (\ref{eq35}), (\ref{eq37}) and the Gauss-Codazzi equation for the screen distribution,
\begin{eqnarray}
\eta(\bar{R}(X,Y)Z)&=&(D_{X}^{g}C)(Y,Z)-(D_{Y}^{g}C)(X,Z)\cr
     & & \varphi_{g}(X)C(Y,Z)-\varphi_{g}(Y)C(X,Z).
\end{eqnarray}

Taking into account (\ref{eq40}) and (\ref{eq38}), we get 

\begin{prop}
\label{ricci}
The Ricci curvature of $D$ is given by 
\begin{eqnarray}
\label{eq44}
Ric^{D}(X,Y)& =& Ric^{g}(X,Y) -2d\theta_{g}(X,Y)+(1-n)(D^{g}_{X}\theta_{g})(Y)\cr & & +(n-1)\theta_{g}(X)\theta_{g}(Y)+ (1-n)g(X,Y)|\theta_{g}^{\sharp_{g}}|^{2}_{g}\cr & &-g(X,Y)\delta^{g}\theta_{g}+ \left([(D^{g}_{X}S)(\xi,Y)-(D^{g}_{\xi}S)(X,Y)]\right.\cr & & \cr& &\left.  +g(X,Y)S(\xi,\theta_{g}^{\sharp_{g}})-S(\xi,X)S(\xi,Y)+\varphi_{g}(X)S(\xi,Y)   \right)
\end{eqnarray}
\end{prop}

\begin{prop}
\label{scal}
Let D be a Weyl structure on $(M,c,S(TM))$, then, for $g \in c$,

\begin{eqnarray}
\label{eq42}
Scal^{D}_{g}&=&scal^{g} - (n-1)^{2}|\theta_{g}^{\sharp_{g}}|^{2}_{g} + (1-2n)\delta^{g}\theta_{g} + (n-1)\varphi_{g}(\theta_{g}^{\sharp_{g}}) \cr & &+ div^{g}i_{\xi}S-tr{g}(D^{g}_{\xi}S)+nS(\xi,\theta_{g}^{\sharp_{g}}) \cr  & &- |(i_{\xi}S)^{\sharp_{g}}|_{g}^{2} +  g(\varphi_{g}^{\sharp_{g}},(i_{\xi}S)^{\sharp_{g}} )    
\end{eqnarray}
\end{prop}

{\bf{Proof}~} We have
$$Scal^{D}_{g}~=~ g^{[\alpha \beta]}Ric^{D}(X_{\alpha},X_{\beta})$$
where $(X_{\alpha})_{\alpha}$ is a quasiorthonormal frame field on $M$ adapted to the decomposition (\ref{eq3}). Then using the above Ricci formula leads to, 

\begin{eqnarray*}
Ric^{D}(X_{\alpha},X_{\beta})&=& Ric^{g}(X_{\alpha},X_{\beta}) -2d\theta_{g}(X_{\alpha},X_{\beta}) + (1-n)(D^{g}_{X_{\alpha}}\theta_{g})(X_{\beta})\cr & & + (n-1)\theta_{g}(X_{\alpha})\theta_{g}(X_{\beta})+ (1-n)g_{\alpha \beta}|\theta_{g}^{\sharp_{g}}|_{g}^{2}-g_{\alpha \beta}\delta^{g}\theta_{g} \cr & &+\left( [(D^{g}_{X_{\alpha}}S)(\xi,X_{\beta})-(D^{g}_{\xi}S)(X_{\alpha},X_{\beta})] + g_{\alpha \beta}S(\xi,\theta_{g}^{\sharp_{g}})\right. \cr & & \left.-S(\xi,X_{\alpha})S(\xi,X_{\beta})+ \varphi_{g}(X_{\alpha})S(\xi,X_{\beta}) \right).
\end{eqnarray*}

with~ $\delta^{g}\theta_{g}:=div^{g}\theta_{g}^{\sharp_{g}}$. Contracting with $g^{[\alpha \beta]}$ and a straighforward computation give relation (\ref{eq42}).

\section{ Einstein-Weyl screen structures}
\label{einsweyl}
Note that as $D$ is not a metric connection on $M$, its Ricci curvature is not necessarily symmetric. The quadruplet $(M,c,S(TM),D)$ defines an Einstein-Weyl screen structure if $D$ is a Weyl screen structure on $(M,c,S(TM))$ and the symmetrised Ricci tensor of $D$ is proportional to $g$ pointwise. Equivalently, there exist a function $\wedge \in C^{\infty}(M)$ such that
\begin{equation}
\label{eq43} 
Ric^{D}(X,Y)~+~Ric^{D}(Y,X) ~=~ \wedge g(X,Y),
\end{equation}
for all tangent vectors $X$, $Y$ $\in TM$. The function ~$\wedge$  (depends on $g\in c$ and ) is called the Einstein-Weyl function of the structure with respect to $g$. 

By (\ref{eq44}) one has 
\begin{eqnarray}
\label{eq45}
Ric^{D}(X,Y)+R^{D}ic(Y,X)&=&Ric^{g}(X,Y)+R^{g}ic(Y,X) +{\mathcal{D}}(\theta_{g})(X,Y)\cr & & + 2g(X,Y)\left\{(1-n)|\theta_{g}^{\sharp_{g}}|_{g}^{2}-\delta^{g}\theta_{g}\right. \cr & & \left.  +S(\xi,\theta_{g}^{\sharp_{g}})  \right\}
\end{eqnarray}

where

\begin{eqnarray}
\label{eq46}
{\mathcal{D}}(\theta_{g})(X,Y)&=&(1-n)\left[(D_{X}^{g}\theta_{g})(Y)+(D_{Y}^{g}\theta_{g})(X)-2\theta_{g}(X)\theta_{g}(Y)\right]\cr & & + \left[(D_{X}^{g}S)(\xi,Y)+(D_{Y}^{g}S)(\xi,X)\right]\cr & & +\left[\varphi_{g}(X)S(\xi,Y)+\varphi_{g}(Y)S(\xi,X)\right]\cr & &-2\left[(D_{\xi}^{g}S)(X,Y) + S(\xi,X)S(\xi,Y)\right]
\end{eqnarray}
Also, on the symmetry of $Ric^{g}$ note that
\begin{equation}
\label{47}
Ric^{g}(X,Y)-Ric^{g}(Y,X)~=~2d\varphi_{g}(X,Y).
\end{equation}
for all tangent vectors $X$, $Y$ in $TM$.
Then, it follows (\ref{eq46}) and (\ref{47})

\begin{prop}
\label{einwecond}
The quadruplet $(M,c,S(TM),D)$ defines a Einstein-Weyl screen  structure if and only if $D$ is defined by (\ref{eq28}) for all $g\in c$ and the Ricci curvature of $g$ satisfies
\begin{equation}
\label{eq48bis}
Ric^{g}~=~d\varphi_{g} - \frac{1}{2}{\mathcal{D}}(\theta_{g}) + \bar{\wedge}g,
\end{equation}
where $\bar{\wedge}$ is related to $\wedge$ in (\ref{eq43}) by
\begin{equation}
\label{eq49bis}
\bar{\wedge} = \frac{1}{2}\wedge - \left[(1-n)|\theta_{g}^{\sharp_{g}}|_{g}^{2}- \delta^{g}\theta_{g} + S(\xi,\theta_{g}^{\sharp_{g}})    \right]
\end{equation}
with ${\mathcal{D}}(\theta_{g})$~ given by (\ref{eq46}).
\end{prop}

\section{Totally umbilical screen foliation.}
\label{totumb}
The screen distribution  $S(TM)$ is said to be totally umbilical if there exist a function $\lambda \in C^{\infty}(M)$ such that 
\begin{equation}
\label{eq50}
C(X,PY)~=~\lambda g(X,Y),
\end{equation}
for all tangent vectors $X$, $Y$ in $TM$. Then,  (\ref{eq29})  becomes
\begin{equation}
\label{eq51}
S(X,Y)=\lambda g(X,Y) + \eta(X)\theta_{g}(Y) ,
\end{equation}
for $(X,Y) \in \Gamma(TM)\times \Gamma(S(TM))$.

In particular,

\begin{equation}
\label{eq52}
S(\xi,X)=:S(X,\xi)~=~\theta_{g}(X)
\end{equation}
for all $X$ in $\Gamma(TM)$.

\begin{lem}
\label{dgzs}
For $g\in c$ and for all  tangent vectors $X$, $Y$ and $Z$ in $\Gamma(TM)$,
\begin{eqnarray}
\label{eq53}
\left(D^{g}_{Z}S\right)(X,Y)&=&(Z\cdot \lambda)g(X,Y)+\left[\theta_{g}(X)(D^{g}_{Z}\eta)(Y)+\theta_{g}(Y)(D^{g}_{Z}\eta)(X) \right]\cr &&  \cr & &+ \left[ \eta(X)(D^{g}_{Z}\theta_{g})(Y) +\eta(Y)(D^{g}_{Z}\theta_{g})(X) \right]
\end{eqnarray}
\end{lem}

\noindent
{\bf{Proof.~}}Let $(X,Y) \in \Gamma(TM)\times \Gamma(S(TM))$, it is immediate using (\ref{eq51}) that, for $Z\in \Gamma(TM)$,  
\begin{equation}
\label{eq54}
\left(D^{g}_{Z}S\right)(X,Y)~=~(Z\cdot \lambda)g(X,Y)+(D^{g}_{Z}\eta)(X)\theta_{g}(Y)+(D^{g}_{Z}\theta_{g})(Y)\eta(X).
\end{equation}
Now for $(X,Y) \in \Gamma(TM)\times \Gamma(TM)$ observe that

$$\left(D^{g}_{Z}S\right)(X,Y)~=~\left(D^{g}_{Z}S\right)(X,PY)+\eta(Y)\left(D^{g}_{Z}S\right)(\xi,PX)$$
and then, using (\ref{eq54}) and the fact that $\theta_{g}$ and $\eta$ are horizontal and vertical respectively, lead to relation (\ref{eq53}).

In particular, for all tangent vectors $X$, $Y$ in $TM$,  

\begin{equation}
\label{eq55}
\left(D^{g}_{X}S\right)(\xi,Y)~=~(D^{g}_{X}\theta_{g})(Y)-\varphi_{g}(X)\theta_{g}(Y).
\end{equation}
and 
\begin{equation}
\label{eq56}
\left(D^{g}_{\xi}S\right)(X,Y)~=~(\xi \cdot\lambda)g(X,Y)+ (D^{g}_{\xi}\theta_{g})(X)\eta(Y)+ (D^{g}_{\xi}\theta_{g})(Y)\eta(X),
\end{equation}
which arises from (\ref{eq53}) and the fact that $\eta$ is parallel along the $\xi-$orbits.

we also have the following fact.
\begin{prop}
\label{facts}
Assume that $(M,c,S(TM),D)$ is an Einstein-Weyl screen structure with totally umbilical $S(TM)$, then
\begin{itemize}
\item[(i)] 
\begin{equation}
\label{eq57}
(D^{g}_{\xi}\theta_{g})(X)~=~0, ~\forall ~ X\in \Gamma(TM).
\end{equation}
and 
\begin{equation}
\label{eq58}
(D^{g}_{\xi}S)(X,Y)=(\xi \cdot\lambda)g(X,Y)
\end{equation}
for all tangent vectors ~$X$,$Y$ in  $\Gamma(TM)$.
\item[(ii)]
\begin{eqnarray}
\label{eq59}
Ric^{D}(X,Y)&=&Ric^{g}(X,Y)-2d\theta_{g}(X,Y)+(2-n)(D_{X}^{g}\theta_{g})(Y)\cr & &+ (n-2)\theta_{g}(X)\theta_{g}(Y)+(2-n)|\theta_{g}^{\sharp_{g}}|_{g}^{2}g(X,Y)\cr & &-g(X,Y)\delta^{g}\theta_{g}-(\xi \cdot \lambda)g(X,Y).
\end{eqnarray}

\item[(iii)]
\begin{eqnarray}
\label{eq60}
Scal^{D}_{g}&=&scal^{g}+(2-n)(n-1)|\theta_{g}^{\sharp_{g}}|_{g}^{2} +2(1-n)\delta^{g}\theta_{g} + n\varphi_{g}(\theta_{g}^{\sharp_{g}})\cr & &-n(\xi \cdot \lambda).
\end{eqnarray}
\end{itemize}
where $\lambda$ is given by (\ref{eq50}).
\end{prop}

\noindent
{\bf{Proof.~}}Note that ~$Ric^{g}(\xi,Y)~=~Ric^{g}(Y,\xi)~=~0~$ and  $~2d\varphi_{g}(\xi,Y)~= Ric^{g}(\xi,Y)-Ric^{g}(Y,\xi) =~0~$. Then, as the structure is Einstein-Weyl, by (\ref{eq48bis}), we have ~${\mathcal{D}}(\theta_{g})(\xi,X)=0~$ for all tangent vector ~$X$~ in ~$\Gamma(TM)$. Thus, (\ref{eq57}) follows (\ref{eq56}) setting $Y=\xi$, and sustitution in (\ref{eq46}).
Thereafter, (\ref{eq56}) reduces to (\ref{eq58}) and $(i)$~is proved.
Now, (\ref{eq59}) and (\ref{eq60}) are just rewriting of (\ref{eq44}) and (\ref{eq42})~respectively, taking into account $(i)$ and (\ref{eq51}) and the proof is complete.                     

{\bf{Note.~}} All metric $g\in c=[g_{0}]_{0}$ will be called the 
\emph{trivial extension} of its restriction $g'\in c'=[g_{0}']$ on the horizontal.
\begin{lem}
\label{ric=ric'}
If $(M,g)$ is totally geodesic in flat $(\ov{M},\ov{g})$ then for all horizontal vector fields $X$ and $Y$, one has 
\begin{equation}
\label{eq61}
Ric^{g}(X,Y)~=~Ric^{g'}(X,Y),
\end{equation}
where $g'$ is the restriction of $g$ on the horizontal.
\end{lem}

{\bf{Proof.~}}For horizontal vector fields $X$~and ~$Y$, one has
\begin{eqnarray*}
Ric^{g}(X,Y)&=&g^{[\alpha \beta]}\tilde{g}(R^{g}(X,X_{\alpha})Y,X_{\beta})\cr
 & &\cr
& =&g^{ij}g(R^{g}(X,X_{i})Y,X_{j})+\tilde{g}(R^{g}(X,\xi)Y,\xi)
\end{eqnarray*}
On the other hand, for horizontal $X$,$Y$~and~$Z$ one has 

\begin{eqnarray*}
R^{g}(X,Y)Z &=& \stackrel{\star}{R}(X,Y)Z+
\left\{ [(\na^{g'}_{Y}C)(X,Z)-(\na^{g'}_{X}C)(Y,Z)] \right.   
\cr & & \cr & & \left. + [C(X,Z)\varphi_{g}(Y)-C(Y,Z)\varphi_{g}(X)]\right\}\xi
\end{eqnarray*}
where $\stackrel{\star}{R}$ denotes the curvature tensor of the induced Levi-Civita connction $\na^{g'}$ on the horizontal. 
Hence, 
\begin{eqnarray*}
Ric^{g}(X,Y)&=&g^{ij}g(\stackrel{\star}{R}(X,X_{i})Y,X_{j})
+\tilde{g}(R^{g}(X,\xi)Y,\xi)\cr & & \cr
& =&g^{'ij}g'(\stackrel{\star}{R}(X,X_{i})Y,X_{j})+\tilde{g}(R^{g}(X,\xi)Y,\xi)\cr  & & \cr
&= &Ric^{g'}(X,Y)+\tilde{g}(R^{g}(X,\xi)Y,\xi)
\end{eqnarray*}
Finally, as (M,g) is totally geodesic in $\ov{M}$ which is flat, we have~\cite[page~97 ]{DB}
$$R^{g}(X,\xi)Y=\ov{R}(X,\xi)Y=0.$$
Thus,

\begin{equation}
\label{eq62}
Ric^{g}|_{hor}=Ric^{g'}.
\end{equation}

\begin{rem}
\label{scal=scal'}
Under hypothesis of lemma~{\ref{ric=ric'}}, since $Ric^{g}(\xi,X)=Ric^{g}(X,\xi)$, it follows (\ref{eq62}) that on leaves of the integrable screen distribution, one has 
\begin{equation}
\label{eq63}
scal^{g}|_{M'}~=~scal^{g'}
\end{equation}
where $M'$ is any leaf of $S(TM)$.
\end{rem}

{\bf{Note.~}}For a Weyl screen structure $D$ relative to $(M,c,S(TM))$, as for any $g\in c $, the associate $1-$form $\theta_{g}$ is horizontal,we will indistinctly note by $\theta_{g}$ its restriction on the horizontal. Thus, for horizontal vectors $X$, $Y$, we have
\begin{equation}
\label{eq64}
\left(D^{g}_{X}\theta_{g} \right)(Y)=X\cdot \theta_{g}(Y)-\theta_{g}(\na^{g'}_{X}Y + C(X,Y)\xi)=\left(\na^{g'}_{X} \theta_{g}\right)(Y).
\end{equation}
Now, we state the following.
\begin{theo}
\label{theo3}
Let $(M,c,S(TM),D)$ be an Einstein-Weyl screen structure quadruplet in the Lorentzian space $\mathbb{R}^{n+2}_{1}$ and $D'$ the (Riemannian) induced Weyl structure by $D$ on the conformal structure $(M',c')$~where $M'$ is a leaf of the totally umbilical integrable screen distribution $S(TM)$and $c'=c|_{M'}$. Then,
\begin{itemize}
\item[(a)] $D'$ is a (Riemann) Einstein-Weyl structure relative to $(M',c')$. Furthermore, the Einstein-Weyl functions $\wedge$~ and $\wedge'$ relative to $g \in c$ and ~$g'=g|_{M'}\in c'$ ~respectively, are related along $M'$ by
\begin{equation}
\label{eq65}
\frac{1}{2}(\wedge-\wedge')~=~\varphi_{g}(\theta_{g}^{\sharp_{g}}) + 2(\xi \cdot \lambda).
\end{equation}

\item[(b)] If the screen foliation is compact and the Cotton-York tensor \cite{Iva} of $D'$ vanishes identically, then the Weyl screen structure $D$ relative to $(M,c,S(TM))$ is closed.

\item[(c)] Along compact leaves of $S(TM)$, the trivial extension $g$ to $(M,c)$ of the Gauduchon metric~\cite{Gaud}~associated to $(M',c',D')$ satisfies
\begin{itemize}
\item[(i)]
\begin{equation}
\label{eq66}
Scal^{g}-(n+2)|\theta_{g}^{\sharp_{g}}|_{g}^{2}~=~G,
\end{equation}

\item[(ii)]
\begin{equation}
\label{eq67}
Scal^{D}_{g}+n(n-4)|\theta_{g}^{\sharp_{g}}|_{g}^{2}-(3-2n)\varphi_{g}(\theta_{g}^{\sharp_{g}})+n(\xi \cdot \lambda)~=~G
\end{equation}
where $G$ is the Gauduchon constant~\cite{Gaud}.
\end{itemize}
\end{itemize}
\end{theo}

\noindent
{\bf{Proof.~}} Let $X$, $Y$ be horizontal vector fields. By use of (~\ref{eq48bis}),~lemma~\ref{ric=ric'}, (\ref{eq55}) and (\ref{eq58}), we have
\begin{eqnarray*}
Ric^{g'}(X,Y) & =& Ric^{g}(X,Y)\cr
              & =& d\varphi_{g}(X,Y)-\frac{1}{2}{\mathcal{D}}(\theta_{g})(X,Y)\cr & &+\left[\frac{1}{2}\wedge -[(2-n)|\theta_{g}^{\sharp_{g}}|_{g}^{2}-\delta^{g}\theta_{g}-2(\xi \cdot \lambda) ]\right]g(X,Y).
\end{eqnarray*}
where $\wedge$ is the Einstein-Weyl function with  respect to $g \in c$. Hence, from (\ref{eq64}) we have 

\begin{eqnarray*}
Ric^{g'}(X,Y)&=& d\varphi_{g}(X,Y)-\frac{1}{2}{\mathcal{D'}}(\theta_{g})(X,Y)\cr & &+\left[\frac{1}{2}\wedge -[(2-n)|\theta_{g}^{\sharp_{g}}|_{g}^{2}-\delta^{g}\theta_{g}-2(\xi \cdot \lambda) ]\right]g'(X,Y).
\end{eqnarray*}
with 

\begin{equation}
\label{eq68}
{\mathcal{D'}}(\theta_{g})(X,Y)= (2-n)\left[(\na_{X}^{g'}\theta_{g})(Y)+ (\na_{Y}^{g'}\theta_{g})(X)-2\theta(X)\theta(Y)\right]
\end{equation}

The symmetry of the $(0,2)-$tensors $Ric^{g'}$, ${\mathcal{D'}}(\theta_{g})$  and $g'$ leads to $d\varphi_{g}(X,Y) = 0$ and 
\begin{eqnarray}
\label{eq69}
Ric^{g'}(X,Y)&=&-\frac{1}{2}{\mathcal{D'}}(\theta_{g})(X,Y)+\cr & &
+\left[\frac{1}{2}\wedge' -[(2-n)|\theta_{g}^{\sharp_{g}}|_{g'}^{2}-\delta^{g'}\theta_{g'}]\right]g'(X,Y).
\end{eqnarray}
with $$\wedge'= \wedge -2\left[ \varphi_{g}(\theta_{g}^{\sharp_{g}}) + 2(\xi \cdot \lambda)  \right].$$
It follows~ (\ref{eq69})~ that $(M',c',D')$ is an Einstein-Weyl structure on the Riemannian leaf $M'$~\cite{Pedswan} with Einstein-Weyl function $\wedge'$ relative to $g'$ as given in (\ref{eq65}).

Now, let $g\in c$ denote the trivial extension of the standard metric of $(M',c',D')$ and $\theta_{g}$ the associated $1-$form. We show that $d\theta_{g}=0$. Suppose $M'$ is a compact leaf of $S(TM)$ and that the Cotton-York tensor of $D'$ vanishes identically. Then, we know by Ianov result in \cite{Iva} that $\na^{g'}\theta_{g}=0$ where $\na^{g'}$ is the Levi-Civita connection of $g'=g_|{M'}$ the standard metric of  $(M',c',D')$.
Then using (\ref{eq64}), we have for horizontal vector fields $X$~and~$Y$,
$$\left(D^{g}_{X}\theta_{g} \right)(Y)=\left(\na^{g'}_{X} \theta_{g}\right)(Y)~=~0.$$ Finally, using (\ref{eq57}) and the fact that $\theta_{g}$ is horizontal, we deduce that $\theta_{g}$ is parallel with respect to $D^{g}$, that is $D^{g}=0$. Hence, $d\theta_{g}=0$ and $\theta_{g}$is closed and $(b)$ is proved.

Note that we have\cite{Gaud} on $(M',c',D',)$
\begin{equation}
\label{eq70}
Scal^{D'}_{g'}=Scal^{g'}+2(n-1)\delta^{g'}\theta_{g'}-(n-1)(n-2)|\theta_{g'}^{\sharp_{g'}}|_{g'}^{2}.
\end{equation}
Also,  the following relation defines the Gauduchon's constant $G$:

\begin{equation}
\label{eq72}
Scal^{D'}_{g'}+n(n-4)|\theta_{g'}^{\sharp_{g'}}|_{g'}^{2}~=~G
\end{equation}

Then, $(i)$ in $(c)$ is a simple consequence of remark~\ref{scal=scal'}. On the other hand, using (\ref{eq60}), remark~\ref{scal=scal'} and (\ref{eq70}), one has along $M'$ 
\begin{equation}
\label{eq73}
Scal^{D}_{g}|_{M'}=Scal^{D'}_{g'}-4(n-1)\delta^{g'}\theta_{g}+(3-2n)\varphi_{g}(\theta_{g}^{\sharp_{g}})-n(\xi \cdot \lambda).
\end{equation}
So, 
\begin{eqnarray}
\label{eq74}
Scal^{D}_{g}|_{M'}+ n(n-4)|\theta_{g}^{\sharp_{g}}|_{g}^{2}+ 4(n-1)\delta^{g'}\theta_{g}\cr -(3-2n)\varphi_{g}(\theta_{g}^{\sharp_{g}})+ n(\xi \cdot \lambda)=      Scal^{D'}_{g'} +n(n-4)|\theta_{g'}^{\sharp_{g'}}|_{g'}^{2}.
\end{eqnarray}
Then, $(ii)$ follows (\ref{eq74}) and (\ref{eq72}) and the proof is complete.

\section{Lightlike real hypersurfaces of Kaehler manifolds}
\label{kahler}
Let $(\ov{M},\ov{g}_{0},\ov{J})$ be a real $2m-$dimensional $(m>1)$ indefinite almost hermitian manifold, where $\ov{g}_{0}$ is a pseudo-Riemannian metric of index $q=2\nu$, $0<\nu <m$. Let $(M,c=[g_{0}]_{0})$ be a lightlike hypersurface of $\ov{M}$ endowed  with the conformal structure $c=[g_{0}]_{0}$~where $g_{0}$ is the degenerate induced metric on $M$ by $\ov{g}_{0}$. As the ambiant manifold $\ov{M}$ has an additional structure $\ov{J}$, it is possible to construct a particular screen distribution on $M$ such that $\ov{J}(TM^{\perp})\oplus \ov{J}(tr(TM))$ be a vector subbundle of $S(TM)$ of rank 2 ~\cite{DB}. More precisely, there exists a nondegenerate almost complex distribution $D_{0}$ with respect to $\ov{J}$ such that 
\begin{equation}
\label{eq75}
S(TM)~=~\left(\ov{J}(TM^{\perp})\oplus \ov{J}(tr(TM))  \right)\perp D_{0}.
\end{equation}
Then, the tangent bundle $TM$ splits as follows
\begin{equation}
\label{eq76}
TM~=~\left(\ov{J}(TM^{\perp})\oplus \ov{J}(tr(TM))  \right)\perp D_{0} \perp TM^{\perp}.
\end{equation}

Consider 

\begin{equation}
\label{eq77}
\Delta~=~\left(\ov{J}(TM^{\perp})\perp  TM^{\perp}) \right)\perp D_{0} ~~\subset ~TM,
\end{equation}
and let $\sigma$ and $Q$ be the projection morphisms of $TM$ on $\Delta$ and $\ov{J}(tr(TM))$, respectively. Also, consider the two isotropic vector fielfs $U=-\ov{J}N$~ and ~$V=-\ov{J}\xi$. Now, define on $M$,

\begin{equation}
\label{eq78}
\theta_{0}(X)~=~g_{0}(X,V)
\end{equation}

and for ~$g=e^{-2f}g_{0}\in c$, associate the $1$-form

\begin{equation}
\label{eq79}
\theta_{g}(X)~=~\theta_{0}(X) + df.
\end{equation}
For all tangent vector $X$ in $TM$, one has
\begin{equation}
\label{eq80}
X~=~\sigma X + \theta_{0}(X)U
\end{equation}
and 
\begin{equation}
\label{eq81}
\ov{J}X~=~FX + \theta_{0}(X)N
\end{equation}
where $F$ is $(1,1)-$tensor globally defined on $M$ by $F=\ov{J}\circ \sigma$. It follows that 

\begin{equation}
\label{eq82}
F^{2}X~=~ -X + \theta_{0}(X)U, \hspace{1cm}   \theta_{0}(U)=1,
\end{equation}
that is $(F,\theta_{0},U)$ defines an almost contact structure on $M$\cite{DB}.

\noindent
{\bf{Note.~}} For the remainder of the text, $(\ov{M},\ov{g}_{0},\ov{J})$ is a Kaehler manifold.

We prove the following technical result.
\begin{lem}
\label{techn}
For all tangent vector fields $X$, $Y$ in $TM$,
\begin{itemize}

\item[(i)]$\left(D^{g_{0}}_{X}\theta_{0}\right)(Y)~=~\theta_{0}(Y)\varphi_{g_{0}}(X)-B(X,FY)$;

\item[(ii)]$\left(D^{g_{0}}_{X}F\right)(Y)~=~\theta_{0}(Y)A_{N}X-B(X,Y)U$;

\item[(iii)]$\varphi_{g_{0}}(X)~=~-\theta_{0}\left(D^{g_{0}}_{X}\right)$;

\item[(iv)]$D^{g_{0}}_{X}\theta_{0}^{\sharp_{g_{0}}}~=~F(\sn X)+\varphi_{g_{0}}(X)\theta_{0}^{\sharp_{g_{0}}}$.
\end{itemize}
\end{lem}

\noindent
{\bf{Proof.~}}
Note that $~\eta(\theta_{0}^{\sharp_{g_{0}}})=0$, and $\theta_{0}^{\sharp_{g_{0}}}= V $. Then we have
\begin{eqnarray}
\label{eq83}
\left(D^{g_{0}}_{X}\theta_{0}\right)(Y)&=&(D^{g_{0}}_{X}g_{0})(Y,\theta_{0}^{\sharp_{g_{0}}})~+~g_{0}(Y,D^{g_{0}}_{X}\theta_{0}^{\sharp_{g_{0}}})\cr
& =&B(\theta_{0}^{\sharp_{g_{0}}},X)\eta(Y)-g_{0}(Y,D^{g_{0}}_{X}\ov{J}\xi).
\end{eqnarray} 
Also, as  $\ov{J}$ is parallel with respect to the Levi-Civita connection on $\ov{M}$, we have
$$D^{g_{0}}_{X}\ov{J}\xi = -\ov{J}(\sn X)-\varphi_{g_{0}}(X)\theta_{0}^{\sharp_{g_{0}}}-B(X,\ov{J}\xi)N.$$
Then, contracting by $Y$ with respect $\ov{g}_{0}$ leads to
$$g_{0}(D^{g_{0}}_{X}\ov{J}\xi,Y)= \ov{g}_{0}(D^{g_{0}}_{X}\ov{J}\xi,Y)=B(X,FY)-\varphi_{g_{0}}(X)\theta_{0}(Y)+B(X,V)\eta(Y).$$
Substituting in (\ref{eq83}) gives $(i)$.

Now, we have 
\begin{equation}
\label{eq84}
(D^{g_{0}}_{X}F)(Y)=\ov{\nabla}_{X}(FY)-B(X,FY)N-F(D^{g_{0}}_{X}Y),
\end{equation}
with $\ov{\nabla}$ the Levi-Civita connection of $(\ov{M},\ov{g}_{0})$. As $FY = \ov{J}Y -\theta_{0}(X)N$, one has by a straightforward use of Gauss-Codazzi formulae
$$\ov{\nabla}_{X}(FY)=F(D^{g_{0}}_{X}Y)+B(X,FY)N+\theta_{0}(Y)A_{N}X -B(X,Y)U, $$ 
whose subtitution in (\ref{eq84})~leads to $(ii)$.

From $(i)$, one has,
$$(D^{g_{0}}_{X}\theta_{0})(U) = \varphi_{g_{0}}(X)-B(X,FU)$$
As  $FU=0$, this is equivalent to 
$$\varphi_{g_{0}}(X)=(D^{g_{0}}_{X}\theta_{0})(U)=X\cdot \theta_{0}(U)- \theta_{0}(D^{g_{0}}_{X}U)=- \theta_{0}(D^{g_{0}}_{X}U)$$
Finally, replacing $Y$ by $\xi$ in $(ii)$, we have
\begin{eqnarray*}
0&=&(D^{g_{0}}_{X}F)(\xi) = D^{g_{0}}_{X}(F\xi)-F(D^{g_{0}}_{X}\xi)\cr 
& & \cr
&= &D^{g_{0}}_{X}(\ov{J}\xi)-F(-\sn X + \varphi_{g_{0}}(X)\xi)\cr
& & \cr
&=& -D^{g_{0}}_{X}\theta_{0}^{\sharp_{g_{0}}}+F(\sn X )+ \varphi_{g_{0}}(X)\theta_{0}^{\sharp_{g_{0}}},
\end{eqnarray*}
which gives $(iv)$~and the proof is complete.

\begin{coro}
\label{coro1}
Let $(M,c=[g_{0}]_{0})$ be a conformal structure on the totally geodesic lightlike hypersurface $(M,g_{0})$ of the Kaehler manifold $(\ov{M},\ov{g}_{0}\ov{J})$. Then,
\begin{itemize}
\item[(i)]$\left(D^{g_{0}}_{X}\theta_{0}\right)(Y)~=~\theta_{0}(Y)\varphi_{g_{0}}(X)$
\item[(ii)]$D^{g_{0}}_{X}\theta_{0}^{\sharp_{g_{0}}}=\varphi_{g_{0}}(X)\theta_{0}^{\sharp_{g_{0}}}$.
\end{itemize} 
\end{coro}

\noindent
{\bf{Proof.~}}Follows lemma~\ref{techn}  and theorem~\ref{theo2}.$\Box$ 

Suppose the second fundamental form of the screen distribution given by (\ref{eq75}) is symmetric on $S(TM)$ so that the latter is integrable. Consider on $(M,c,S(TM))$ the Weyl screen structure $D^{\ov{J}}$ defined by (\ref{eq28}) where for each $g\in c$ the associated $1-$form $\theta_{g}$ is given by (\ref{eq79}). For 
$D^{\ov{J}}$ to be closed we prove the following.
\begin{theo}
\label{theo4}
The Weyl structure $D^{\ov{J}}$ relative to $(M,c,S(TM))$ is closed if and only if the $1-$forms $\theta_{0}$ and $\varphi_{g_{0}}$ are proportional.
\end{theo}

\noindent
{\bf{Proof.~}}
For all tangent vector fields $X$, $Y$ in $TM$, we have

$$d\theta_{0} (X,Y)~=~\frac{1}{2}\left[D^{g_{0}}_{X}\theta_{0}(Y)-D^{g_{0}}_{Y}\theta_{0}(X)   \right]$$
and by use of corollary~\ref{coro1}, one has
\begin{eqnarray*}
d\theta_{0} (X,Y)&=&\frac{1}{2}\left[ \varphi_{g_{0}}(X)\theta_{0}(Y)-\varphi_{g_{0}}(Y)\theta_{0}(X)   \right]\cr 
&=& \frac{1}{2}\left[\varphi_{g_{0}} \wedge  \theta_{0} \right] (X,Y).
\end{eqnarray*}
Thus, $\theta_{0}$ is closed if and only $\theta_{0}$ and $\varphi_{g_{0}}$ are proportional.

\vspace{1.5cm}
{\bf{Acknowledgments.}}
The first named author (C. Atindogbe) thanks the Agence Universitaire de la Francophonie (AUF) for support with a one year research grant, along with the Institut Elie Cartan (IECN, UHP-Nancy~I) for  research facilities during the completion of this work.

\end{document}